\RequirePackage{ifpdf}
\ifpdf 
\documentclass[pdftex]{sigma}
\else
\documentclass{sigma}
\fi

\newcommand{\Aut}{\mathop{\sf Aut}\nolimits}

\newcommand{\tr}{\mathop{\sf tr}\nolimits}

\def \Z {{\mathbb Z}}
\def \R {{\mathbb R}}
\def \C {{\mathbb C}}
\def \H {{\mathbb H}}
\def \tr {{\rm Tr}\,}
\def \mat {{\rm Mat}}
\def \ens {{\rm End}}
\def \id {{\rm Id}}
\def \re {{\Re e}}
\def \char {{\rm char}\,}
\def \F {{\mathbb K}}
\def \oF {\overline{\mathbb K}}
\def \Mu {{\bf \mu}}

\def \lc {\left<}
\def \rc {\right>}

\def \bou {\partial \Omega}

\begin{document}

\allowdisplaybreaks

\renewcommand{\PaperNumber}{070}

\FirstPageHeading

\ShortArticleName{Klein Topological Field Theories from Group Representations}

\ArticleName{Klein Topological Field Theories\\ from Group Representations}

\Author{Sergey A. LOKTEV~$^{\dag\ddag}$  and Sergey M. NATANZON~$^{\dag\ddag\S}$}

\AuthorNameForHeading{S.~Loktev and S.~Natanzon}

\Address{$^\dag$~Department of Mathematics, Higher School of Economics,\\
\hphantom{$^\dag$}~7 Vavilova Str., Moscow 117312, Russia}
\EmailD{\href{mailto:S.Loktev@gmail.com}{S.Loktev@gmail.com}, \href{mailto:natanzons@mail.ru}{natanzons@mail.ru}}
\URLaddressD{\url{http://www.hse.ru/en/org/persons/23004570}}
\URLaddressD{\url{http://www.hse.ru/en/org/persons/14026884}}

\Address{$^\ddag$~Institute of Theoretical and Experimental Physics,\\
\hphantom{$^\ddag$}~25 Bolshaya Cheremushkinskaya Str., Moscow 117218, Russia}

\Address{$^\S$~A.N.~Belozersky Institute,  Moscow State University,\\
\hphantom{$^\S$}~Leninskie Gory 1, Bldg.~40, Moscow 119991, Russia}

\ArticleDates{Received December 15, 2010, in f\/inal form July 04, 2011;  Published online July 16, 2011}

\Abstract{We show that any complex (respectively real) representation of f\/inite group
naturally generates a open-closed (respectively Klein) topological
f\/ield theory over complex numbers. We relate the 1-point correlator for
the projective plane in this theory with the Frobenius--Schur indicator on the representation.
We relate any complex simple Klein TFT to a real division ring.}


\Classification{57R56; 20C05}

\pdfbookmark[1]{Introduction}{intr}

\section*{Introduction}\label{s1}

Topological quantum f\/ield theories were introduced by
Schwarz \cite{Sch}, Atiyah \cite{At}, Segal \cite{Se},
and Witten \cite{W}. In this paper we concentrate on
open-closed and Klein topological f\/ield theories.
There are several approaches how to generalize the TQFT framework
for non-compact and non-orientable surfaces related to
dif\/ferent questions of physics and geometry
(see~\cite{Du, FFFS, FRS, LR,MS,N,TT}).
The setting we describe in this paper was proposed in~\cite{AN} and developed in
\cite{AN1,AN2,AN3}. It generalizes the application of TQFT to
theory of Hurwitz numbers discovered in~\cite{D1}.

In Section~\ref{s2} we reformulate def\/initions of closed, open-closed,
and Klein topological f\/ield theories in useful for us form
(for shortness let us omit the word ``quantum'' there).
We recall that categories of
these theories are equivalent to categories of Frobenius pairs,
Cardy--Frobenius algebras and equipped Cardy--Frobenius algebras
respectively (this is proved in~\cite{D2,AN} and widely generalized in  \cite{N1,N2}). Therefore constructions of
a topological f\/ield theories are reduced to constructions of
(equipped) Cardy--Frobenius algebras.

In Section~\ref{s3} we prove that the group algebra and the center of group
algebra of any f\/inite group $G$ form a semi-simple
equipped Cardy--Frobenius algebra over any number f\/ield. We call it regular.
Also we construct a real equipped Cardy--Frobenius algebra from a real division ring
and show that any semi-simple complex equipped Cardy--Frobenius algebra
can be obtained from it by complexif\/ication.
Then we present full description of Regular complex algebras of a~group in these terms.

In Section \ref{s4} we prove that the center of group algebra together with
the intertwining algebra of any representation of~$G$ generates a Cardy--Frobenius algebra
and that this algebra is
equipped if the representation is real.
For representations, that appear
from group actions, we relate this construction with the one proposed in~\cite{AN3}.

\section{Topological f\/ield theories and related algebras}\label{s2}

\subsection{Closed topological f\/ield theories} 
\begin{figure}[tbph]
\centering
\includegraphics[width=6.65cm]{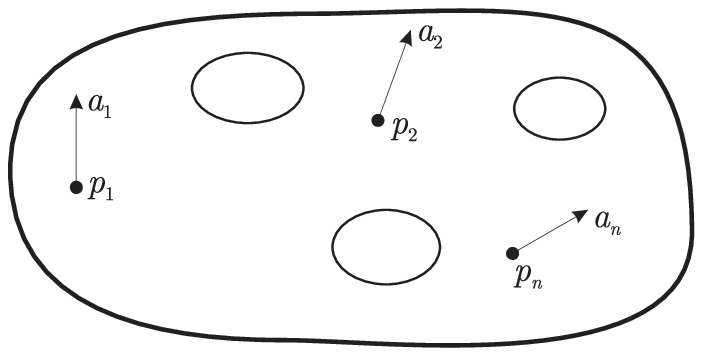}
\caption{}\label{f1}
\end{figure}

The simplest variant
of topological f\/ield theory is closed topological f\/ield theory
(\cite{At, D2}, see~\cite{K} for more references). In this case we consider oriented closed
surfaces without boundary. Also  we f\/ix a f\/inite-dimensional vector space
$A$ over a f\/ield $\mathbb{K}$ with basis $\alpha_1, \dots, \alpha_N$ and associate a~number
$\lc a_1,a_2,\dots ,a_n \rc_{\Omega}$ to each system of vectors
$a_1,a_2,\dots ,a_n\in A$ situated at a set of points
$p_1,p_2,\dots ,p_n$ on a surface $\Omega$ (Fig.~\ref{f1}).

We assume that the numbers $\lc a_1,a_2,\dots ,a_n \rc_{\Omega}$
are invariant with respect to
any homeo\-morphisms of surfaces with marked points. Moreover, we postulate
that the system $\{ \langle a_1,a_2,\dots$, $a_n\rangle_{\Omega}\}$ consists of
multilinear forms and satisf\/ies a non-degeneracy axiom and cut axioms.

The \textit{\textit{non-degeneracy axiom}} says that the matrix
$\left(\lc \alpha_i,\alpha_j \rc_{S^2}\right)_{1\le i,j \le N}$
is non-degenerate.
By $F^{\alpha_i,\alpha_j}_A$ denote
the inverse matrix.

The \textit{cut axioms}  describes evolution of
$\lc a_1,a_2,\dots ,a_n \rc_{\Omega}$ by cutting and collapsing along  contours
$\gamma\subset\Omega$. Indeed, there are two cut axioms related to dif\/ferent
topological types of contours.

\begin{figure}[tbph]
\centering
\includegraphics[width=8.55cm]{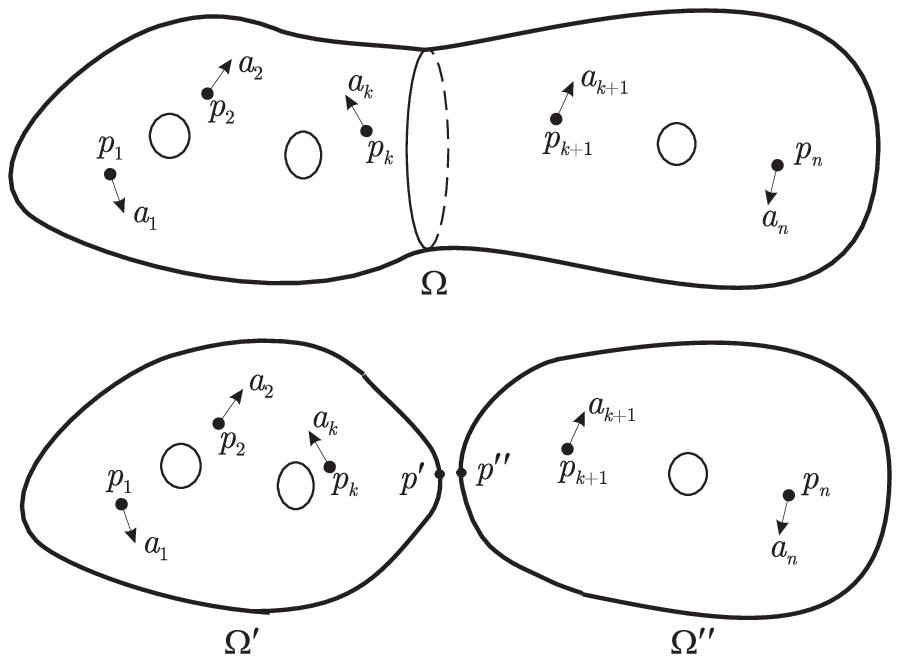}
\caption{}\label{f2}
\end{figure}

If $\gamma$ decomposes $\Omega$ into  $\Omega'$ and
$\Omega''$ (Fig.~\ref{f2})

\noindent
then
\begin{gather*}
\lc a_1,a_2,\dots ,a_n \rc_{\Omega} = \sum_{i,j}
\lc a_1,a_2,\dots ,a_k,\alpha_i \rc_{\Omega'}F^{\alpha_i,\alpha_j}_A
\lc \alpha_j,a_{k+1}, a_{k+2},\dots ,a_n \rc_{\Omega''}.
\end{gather*}

\begin{figure}[tbph]
\centering
\includegraphics[width=9cm]{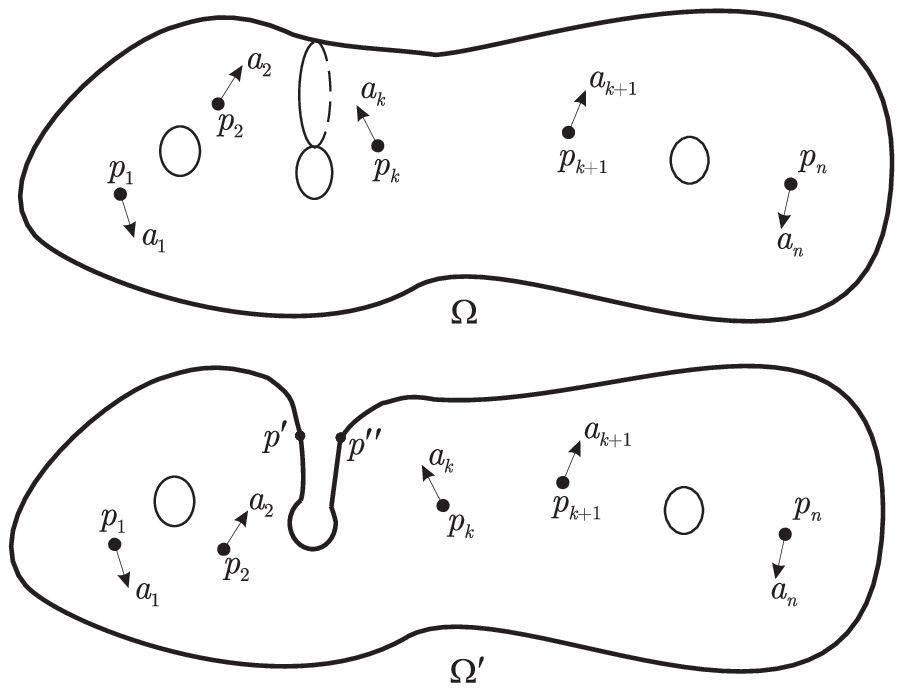}
\caption{}\label{f3}
\end{figure}

If $\gamma$ does not decompose $\Omega$ (Fig.~\ref{f3})

\noindent
then
\begin{gather*}
\lc a_1,a_2,\dots ,a_n \rc_{\Omega} = \sum_{i,j}
\lc a_1,a_2,\dots ,a_n,\alpha_i,\alpha_j \rc_{\Omega'}F^{\alpha_i,\alpha_j}_A.
\end{gather*}

\begin{figure}[tbph]
\centering
\includegraphics[width=10cm]{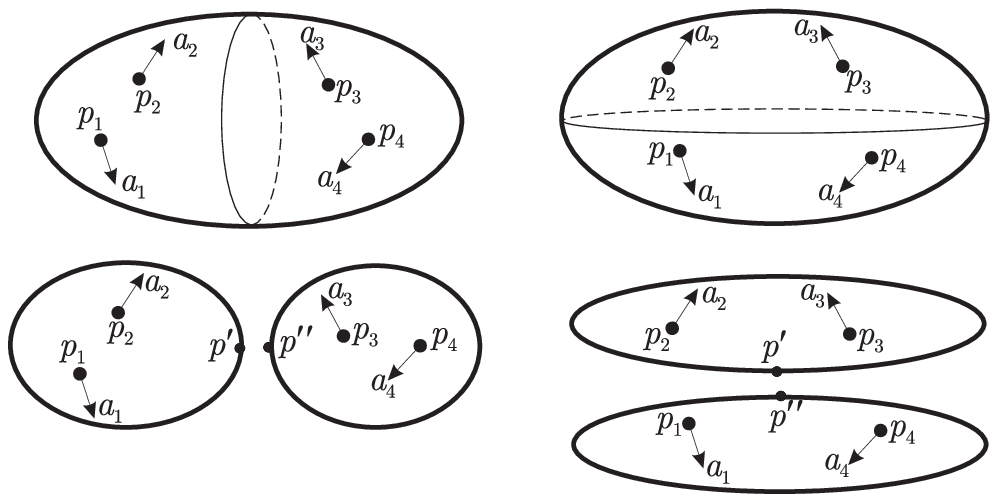}
\caption{}\label{f4}\vspace{-6mm}
\end{figure}

The f\/irst consequence of the topological f\/ield theory axioms is a
\textit{structure of algebra} on~$A$. Namely, the multiplication is def\/ined by
$\lc a_1 a_2,  a_3 \rc_{S^2}= \lc a_1, a_2, a_3 \rc_{S^2}$, so the numbers
$c_{ij}^{k}=\sum_{s} \lc \alpha_i, \alpha_j, \alpha_s \rc_{S^2} F^{\alpha_s,\alpha_k}_A$
are structure constants for this algebra. The cut axiom gives (Fig.~\ref{f4})

\begin{gather*}
\sum_{i,j}
\lc a_1,a_2,\alpha_i \rc_{S^2} F_A^{\alpha_i,\alpha_j} \lc \alpha_j,a_3,a_4 \rc_{S^2}= \lc a_1,a_2,a_3,a_4 \rc_{S^2}\\
\qquad{} = \sum_{i,j}
\lc a_2,a_3,\alpha_i \rc_{S^2}F_A^{\alpha_i,\alpha_j} \lc \alpha_j,a_4,a_1 \rc_{S^2}.
\end{gather*}

Therefore
$\sum_{s,t}c_{ij}^{s}c_{sk}^{t}=\sum_{s,t}c_{jk}^{s}c_{si}^{t}$ and
thus~$A$ is an associative algebra. The unit of the algebra $A$ is given by the vector
$\sum_{i} \lc \alpha_i \rc_{S^2} F^{\alpha_i,\alpha_j}_A\alpha_j$.
The linear form $l(a)=\langle a\rangle_{S^2}$ is a co-unit, also it def\/ines
a non-degenerate invariant bilinear form
$(a_1,a_2)_{A}=l(a_1a_2)= \lc a_1,a_2 \rc_{S^2}$ on $A$. The topological invariance
makes all marked points $p_i$ equivalent and, therefore, $A$~is a~commutative algebra.

Following~\cite{AN2} we say that $(A,l_A)$ is a {\em Frobenius pair} if $A$ is a symmetric
Frobenius algebra~\cite{F}, that is an algebra with
a unit and an invariant non-degenerate symmetric bilinear form, and
the bilinear form is given by $(a_1,a_2) = l_A(a_1\cdot a_2)$.

Thus, $(A,l_A)$ is a commutative Frobenius pair,
moreover, the construction gives a functor $\mathcal{F}$ from the
category of closed topological f\/ield theories to the category of
commutative Frobenius pairs.

\begin{theorem}[\cite{D2}] The functor $\mathcal{F}$ is an equivalence
between categories of closed topological field theories and commutative
Frobenius pairs.
\end{theorem}

The Frobenius structure gives an explicit formula for correlators:
\[
\lc a_1,a_2,\dots ,a_n \rc_{\Omega}=l_A\big(a_1a_2\cdots a_n \left(K_A\right)^g\big),
\] where
$K_A=\sum_{ij}F^{\alpha_i,\alpha_j}_A\alpha_i\alpha_j$ and $g$ is
genus of $\Omega$.

\begin{note}
This def\/inition of closed topological f\/ield theory is equivalent
to the standard one~\cite{At} as follows. Our points $p_1, \dots, p_n$
correspond to circles in a closed 1-dimensional variety,
and the surface corresponds to the cobordism of this variety with
the empty set. The following def\/inition of open-closed topological f\/ield theories is equivalent to proposed in \cite{Laz,MS} in a similar way.
\end{note}

\subsection{Open-closed topological f\/ield theories}

\begin{figure}[tbph]
\centering
\includegraphics[width=4cm]{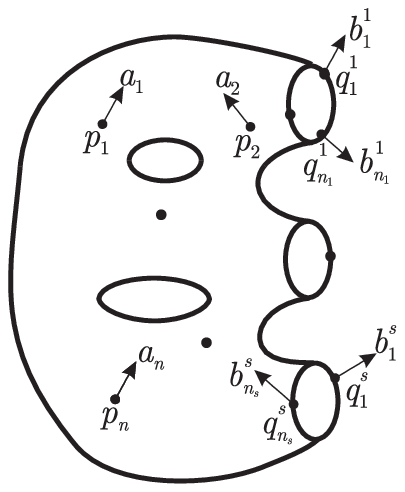}
\caption{}\label{f5}
\end{figure}

A more complicated variant
of topological f\/ield theory is open-closed topological f\/ield theory
\cite{AN,Laz,Moore,MS}. In this case we admit
oriented
compact surfaces $\Omega$  with boundary $\bou$ and some
marked points on $\bou$.
Let us denote the interior marked points and vectors
as before by $p_1,p_2,\dots ,p_n$ and $a_1,a_2,\dots ,a_n\in A$.
Here we introduce a special notation $q^i_j$ for the boundary marked points, where
the upper index $i=1\dots s$ corresponds to a connected component of~$\bou$ (that is a boundary
contour of $\Omega$). The lower index~$j$ is individual for any
boundary contour, it counts the points consequently on the circle,
following the direction determined by the orientation of $\Omega$.
The
vectors $b^i_j$ attached to $q^i_j$
belong to another vector space $B$ over $\mathbb{K}$ with basis $\beta_1, \dots, \beta_M$
(Fig.~\ref{f5}).

Let us denote the corresponding correlator
by  $\lc a_1, \dots, a_n, (b_1^1, \dots, b_{n_1}^1), \dots, (b^s_1, \dots,b^s_{n_s}) \rc_\Omega$
to  keep in mind this picture.
Note that dif\/feomorphisms of $\Omega$ can induce any permutation of $a_i$, but only cyclic
permutations in each group $b_1^i, \dots, b_{n_i}^i$.

We suppose that topological invariance axiom and all axioms of
closed topological f\/ield theory are fulf\/illed for interior marked
points and cut-contours. Thus open-closed topological f\/ield theory
also generates a commutative Frobenius pair $(A, l_A)$. Also we impose
an additional non-degeneracy axiom
and cut axioms related to the boundary.

\textit{The additional non-degeneracy axiom} says that
for any disk  $D$  with two marked boundary points
the matrix $\left( \lc \beta_i,\beta_j\rc_{D} \right)$,
is non-degenerate.
By $F^{\beta_i,\beta_j}_B$  denote the inverse matrix.
It play for ``segment-cuts'' the same ``gluing
role'' that $F^{\alpha_i,\alpha_j}_A$ for ``contour-cuts''.

\begin{figure}[tbph]
\centering
\includegraphics[width=9cm]{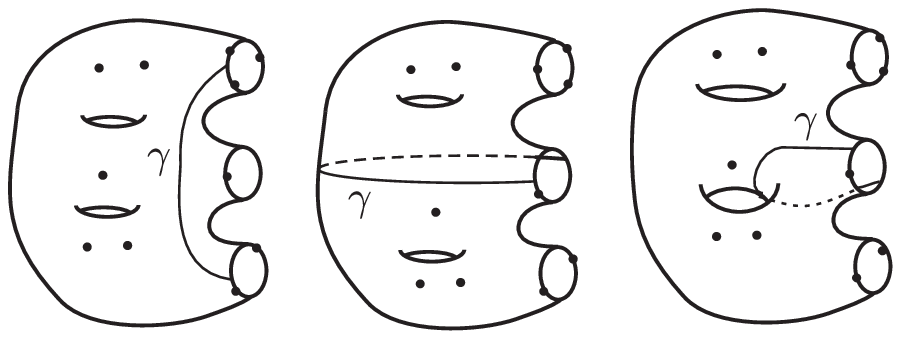}
\caption{}\label{f6}
\end{figure}

In open-closed topological f\/ield theory we consider cuts by
simple segments $[0,1] \to \Omega$ such that the image of
$0$ and $1$ belongs to the boundary.
Then there are three topological types of such cuts (Fig.~\ref{f6}).

\begin{figure}[tbph]
\centering
\includegraphics[width=9cm]{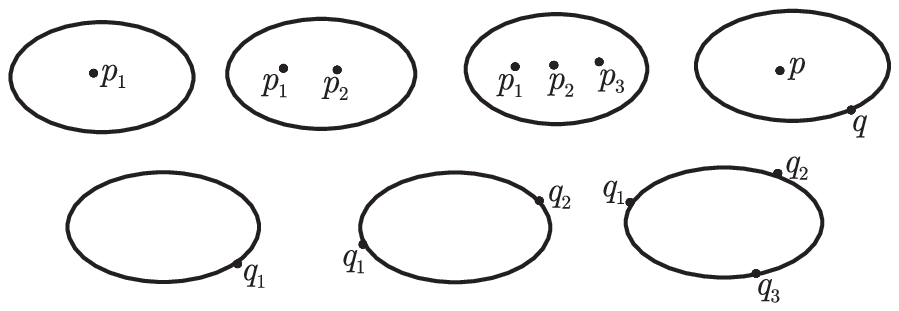}
\caption{}\label{f7}
\end{figure}

Using such cuts one can reduce any marked oriented
surface to elementary marked surfaces from the next list (Fig.~\ref{f7}).
Four surfaces from the upper row are spheres and the remaining three are disks.

\begin{figure}[tbph]
\centering
\includegraphics[width=8cm]{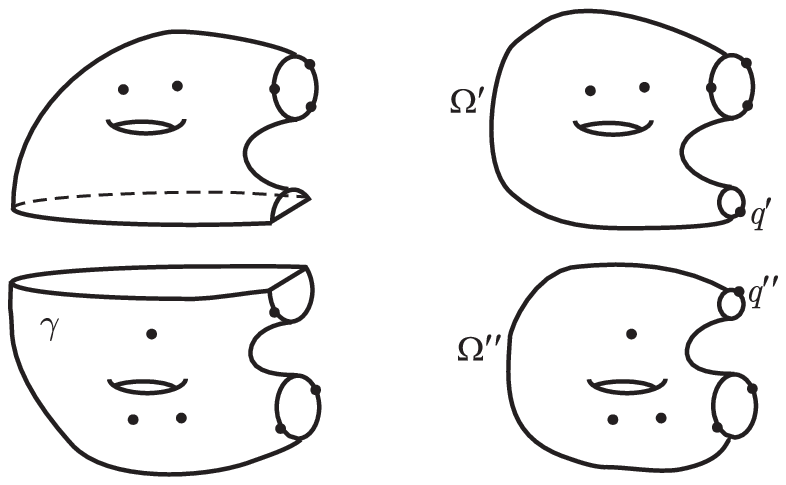}
\caption{}\label{f8}
\end{figure}

Three topological types of segments provide three  new cut axioms.
For example, the axiom for the cut of type 2 (Fig.~\ref{f8}) is

\vspace*{-6mm}

\begin{gather*}
\lc a_1,a_2,\dots ,a_n, (b^1_1,\dots,b^1_{n_1}),  \dots,
(b^s_1,\dots,b^s_{n_s}) \rc_{\Omega}\\
\qquad{} =\sum_{i,j}
\lc a_1,a_2,\dots,a_k, (b^1_1,\dots,b^1_{n_1}), \dots,
(b^t_1,\dots,b^t_{n_t'},\beta_i) \rc_{\Omega'} F^{\beta_i,\beta_j}_B\\
\qquad\quad{}\times
\lc a_{k+1},\dots,a_n, (\beta_j,b^t_{n_t'+1},\dots,b^t_{n_t}), \dots,
(b^s_1,\dots,b^s_{n_s}) \rc_{\Omega''}.
\end{gather*}

The correlators for the  disk $D$ with up to three  boundary points $\lc (b_1) \rc_D$,
$\lc (b_1,b_2) \rc_D$ and $\lc (b_1,b_2,b_3)\rc_D$ give us  a  Frobenius
pair $(B,l_B)$ with structure constants def\/ined in a usual way:
$d_{ij}^{k}=\sum_{s} \left< (\beta_i,\beta_j,\beta_s) \right>_{D} F^{\beta_s,\beta_k}_B$.

\begin{figure}[tbph]
\centering
\includegraphics[width=9cm]{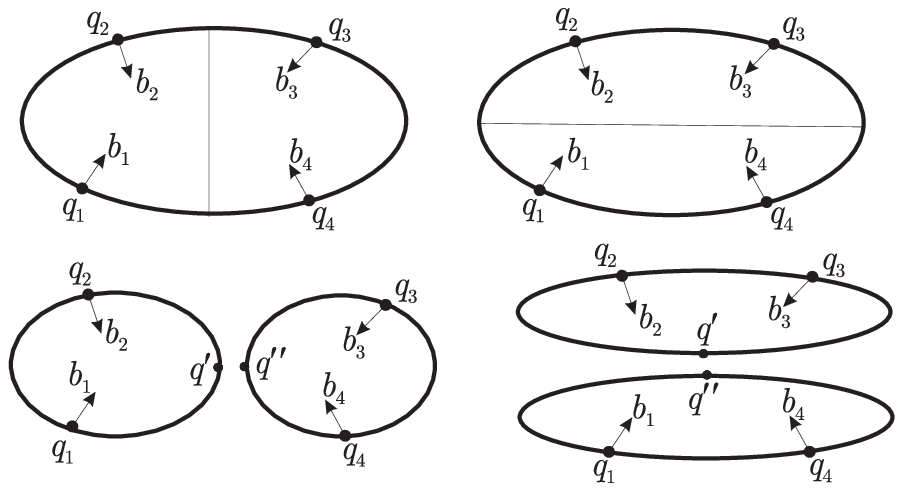}
\caption{}\label{f9}
\end{figure}

The associativity of $B$ follows from the picture below (Fig.~\ref{f9})

Thus $\lc (b_1,b_2,b_3,b_4) \rc_{D}$ is equal both to
$\sum_{i}\lc (b_1,b_2,\beta_i) \rc_{D}F_B^{\beta_i,\beta_j} \lc(\beta_j,b_3,b_4)\rc_{D}$
as well as to
$\sum_{i} \lc (b_2,b_3,\beta_i)\rc_{D}F_B^{\beta_i,\beta_j} \lc (\beta_j,b_4,b_1) \rc_{D}$.
However the algebra $B$ is not commutative in general,
because there is no homeomorphisms of disk that interchanges
$q_1$ with $q_2$ and preserves $q_3$. Nevertheless, calculating corellators for
two boundary points, we have
$l_B(b_1 b_2) = l_B(b_2 b_1)$ that makes the bilinear form on $B$ symmetric.

The correlator $\lc a,(b) \rc_D: A\times B\rightarrow\mathbb{C}$ together with
non-degenerate bilinear forms $\lc a_1,a_2 \rc_{S^2}: A\times
A\rightarrow\mathbb{C}$, $\lc (b_1,b_2) \rc_D: B\times
B\rightarrow\mathbb{C}$ generates two \textit{homomorphisms of vector
spaces} $\phi: A\rightarrow B$ and the dual one $\phi^*: B\rightarrow A$.

Let us deduce some consequences from additional topological axioms.

\begin{proposition}
We have
\begin{enumerate}\itemsep=0pt
\item[$1)$] $\phi$ is an algebra homomorphism,

\item[$2)$] $\phi(A)$ belongs to center of $B$,

\item[$3)$] $(\phi^*(b'),\phi^*(b''))_A=\tr W^B_{b'b''}$, where the
operator $W^B_{b'b''}:B\rightarrow B$ is $W^B_{b'b''}(b)=b'bb''$.
\end{enumerate}
\end{proposition}

The last condition  is known as
the \textit{Cardy condition}, it
appeared in papers  \cite{C1,C2}.

So introduce the notion of a
\textit{Cardy--Frobenius algebra} $((A,l_A),(B,l_B),\phi)$ formed by
\begin{enumerate}\itemsep=0pt
\item[1)] a commutative Frobenius pair $(A,l_A)$;

\item[2)] an arbitrary Frobenius pair $(B,l_B)$;

\item[3)] an algebra homomorphism $\phi: A\rightarrow B$ such that $\phi(A)$
is contained in the center of $B$ and
$(\phi^*(b'),\phi^*(b''))_A=\tr W^B_{b'b''}$.
\end{enumerate}

Thus, we construct a functor $\mathcal{F}$ from the category of
open-closed  topological f\/ield theories to the category of
Cardy--Frobenius algebras $((A,l_A),(B,l_B),\phi)$.

\begin{note}
This conditions are equivalent to  (3.33), (3.44) in \cite{LP}, where $\phi$ is called {\em zipper} and $\phi^*$  is
called {\em cozipper}.
\end{note}

\begin{theorem}[\cite{AN}, see also \cite{LP}]  The functor $\mathcal{F}$
is an equivalence between categories of open-closed topological field
Theories and Cardy--Frobenius algebras.
\end{theorem}

The structure of Cardy--Frobenius algebra provides an explicit formula for correlators:
\begin{gather*}
\lc a_1,a_2,\dots,a_n,\big(b^1_1,\dots,b^1_{n_1}\big),\dots,
(b^s_1,\dots,b^s_{n_s})\rc_{\Omega}\\
\qquad{}= l_B\big( \phi\big(a_1 a_2 \cdots a_n K_A^g\big)  b^1_1 \cdots b^1_{n_1}   V_{K_B}\big(b^2_1 \cdots b^2_{n_2}\big)\cdots
V_{K_B}\big(b^s_1 \dots b^s_{n_s}\big)\big),
\end{gather*}
where the operator
$V_{K_B}:B\rightarrow B$ is given by
$V_{K_B}(b)=F_B^{\beta_i,\beta_j}\beta_i b \beta_j$, and $g$ is genus
of $\Omega$.

\subsection{Klein topological f\/ield theories}

The orientability restriction is indeed avoidable,
the corresponding settings were introduced in~\cite{AN} as Klein topological f\/ield
theory.
It is an extension of open-closed topological f\/ield theory
to arbitrary compact surfaces (possible non-orientable and with
boundary) equipped by a f\/inite set of marked points with local
orientation of their vicinities.

\begin{figure}[tbph]
\centering
\includegraphics[width=3.8cm]{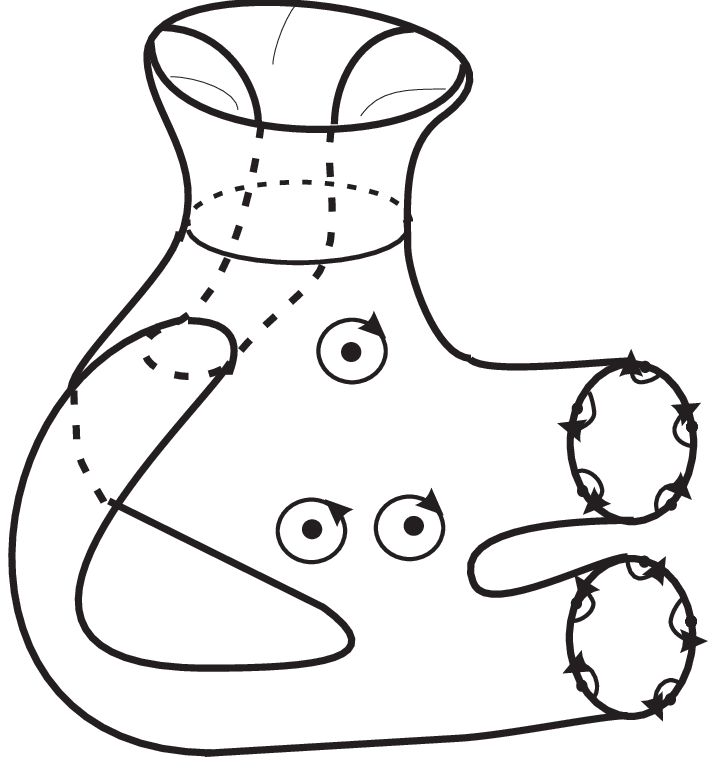}
\caption{}\label{f10}
\end{figure}

In the same way as in open-closed topological f\/ield theory in order to
calculate a correlator we
attach vectors from a space $A$ (respectively $B$)  to interior (resp.\ boundary) marked points
on the surface.

We assume that topological invariance axiom and all axioms of
open-closed topological f\/ield theory are fulf\/illed for cuts that
belong to any orientable part of the surface. Thus Klein topological
f\/ield theory also generates a Cardy--Frobenius algebra
$((A,l_A),(B,l_B),\phi)$.

\begin{figure}[tbph]
\centering
\includegraphics[width=6.65cm]{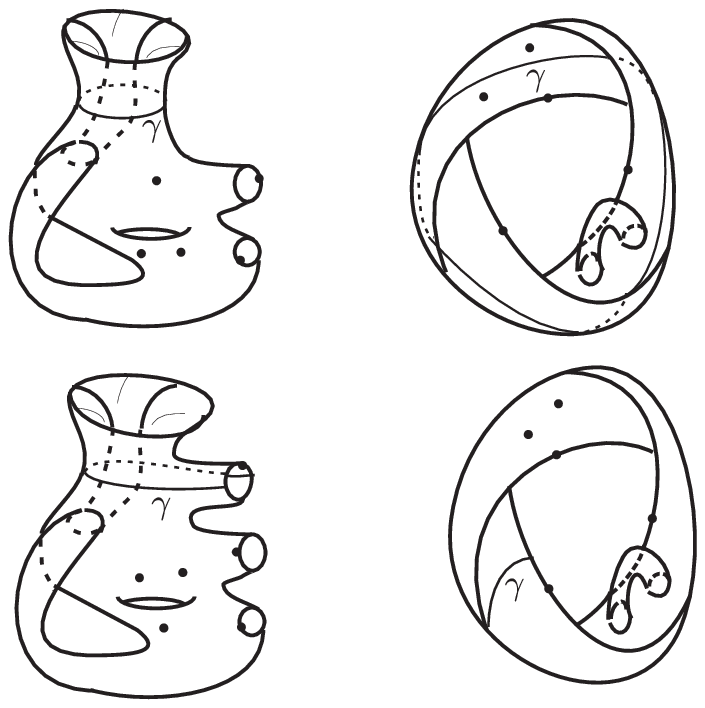}
\caption{}\label{f11}
\end{figure}

Non-orientable surfaces gives 4 new types
of cuts (2 types of cuts by segments and 2 types of cuts by
contours) (Fig.~\ref{f11}).

The full system of cuts provides a way to reduce any marked
non-orientable surface to marked surfaces from list on Fig.~\ref{f7} and
the projective plane with one marked point.

Four new topological type of cuts give four new topological axioms.
Three of them are similar to introduced above, and
the
axiom for the cut of type 2 should be dif\/ferent because
there is only one new contour after the cut.
For this case  introduce an element $U \in A$
such that
\begin{gather*}
 \lc a_1,a_2,\dots,a_n,\big(b^1_1,\dots,b^1_{n_1}\big),\dots,
\big(b^s_1,\dots,b^s_{n_s}\big) \rc_{\Omega}\\
\qquad{} =\lc a_1,a_2,\dots,a_n,U,\big(b^1_1,\dots,b^1_{n_1}\big),\dots,
\big(b^s_1,\dots,b^s_{n_s}\big)\rc_{\Omega'}.
\end{gather*}
Then it follows that $U$ is dual to the 1-point correlator
for the projective plane, namely let $l_U(a)=\lc a \rc_{{\mathbb R \mathbb P}^2}:
A\rightarrow\mathbb{K}$, then $l_A(Ua) = l_U(a)$.

Now suppose that there are linear involutions $\star:A\rightarrow A$ and $\star:B\rightarrow B$,
such that applying~$\star$ inside the correlator gives the same answer as
changing the local orientation around correspon\-ding points.  Sometimes we will write $c^\star=\star(c)$.
Let us  describe the algebraic consequences of these assumptions.

\begin{proposition}[\cite{AN}]
We have
\begin{enumerate}\itemsep=0pt
\item[$1)$] the involution $\star: A \rightarrow A$ is automorphism, the involution
$\star : B \rightarrow B$ is anti-automorphism (that is
$(b_1b_2)^\star =b_2^\star b_1^\star$),

\item[$2)$] $l_A(x^\star)=l_A(x)$, $l_B(x^\star)=l_B(x)$,
$\phi(x^\star)=\phi(x)^\star$,

\item[$3)$] $U^2=F_A^{\alpha_i,\alpha_j}\alpha_i\alpha_j^\star$,
let us denote  it by $K_A^\star$,

\item[$4)$] $\phi(U)=F_B^{\beta_i,\beta_j}\beta_i\beta_j^\star$,
let us denote  it by $K_B^\star$.
\end{enumerate}
\end{proposition}

Summarizing these properties, introduce the notion of
an \textit{equipped Cardy--Frobenius algebra} $((A,l_A),(B,l_B),\phi,U,\star)$ formed by
\begin{enumerate}\itemsep=0pt
\item[1)] a Cardy--Frobenius algebra $((A,l_A),(B,l_B),\phi)$;

\item[2)] an anti-automorphisms $\star :A\rightarrow A$ and $\star : B\rightarrow
B$ such
that $l_A(x^\star)=l_A(x)$, $l_B(x^\star)=l_B(x)$, $\phi(x^\star)=\phi(x)^\star$;

\item[3)] an element $U\in A$ such that $U^2=K_A^\star$ and $\phi(U)=K_B^\star$.
\end{enumerate}

Thus, we constructed a functor $\mathcal{F}$ from the category of
Klein topological f\/ield theory to the category of equipped
Cardy--Frobenius algebras $((A,l_A),(B,l_B),\phi,U,\star)$.

\begin{theorem}[\cite{AN}] The functor $\mathcal{F}$
is an equivalence between categories of Klein topological field
theories and equipped Cardy--Frobenius algebras.
\end{theorem}

The structure of equipped Cardy--Frobenius algebra provides an explicit
formula for
correlators on non-orientable surfaces:
\begin{gather*}
\lc a_1,a_2,\dots,a_n,\big(b^1_1,\dots,b^1_{n_1}\big),\dots,
\big(b^s_1,\dots,b^s_{n_s}\big)\rc_{\Omega}\\
\qquad{}=
l_B\big( \phi\big(a_1 a_2 \cdots a_n U^{2g}\big) b^1_1 \cdots b^1_{n_1} V_{K_B}\big(b^2_1 \cdots b^2_{n_2}\big) \cdots
V_{K_B}\big(b^s_1\dots b^s_{n_s}\big)\big),
\end{gather*}
where $g$ is geometrical genus of $\Omega$, that is, $g=a+1$, if $\Omega$ is the Klein bottle
with $a$ handles and $g=a+\frac{1}{2}$, if $\Omega$ is the projective
plane with $a$ handles.

\section[Regular Cardy-Frobenius algebra of finite group]{Regular Cardy--Frobenius algebra of f\/inite group}\label{s3}

\subsection{Construction of regular algebra}\label{s3.1}

In this section we present a construction that associates
an equipped Cardy--Frobenius algebra and, therefore, a Klein topological f\/ield
theory  to any f\/inite group $G$.

By $|M|$ denote cardinality of a f\/inite set $M$.
Let $\mathbb{K}$ be any f\/ield such that $\char \mathbb{K}$ is not a~divisor of $|G|$.
By $B = \mathbb{K}[G]$ denote the group algebra. It
can be def\/ined as the algebra, formed by linear combinations of
elements of $G$ with the natural multiplication, as well as the
algebra of $\mathbb{K}$-valued functions on $G$ with multiplication def\/ined
by convolution. It has a natural structure of a Frobenius
pair with $l_B(f) = f(1)$. Note that $l_B(f) = {\rm Tr}_{{\mathbb K}[G]} f /|G|$.

The center $A = Z(B)$ with the functional $l_A(f) = f(1)/|G|$ forms
 a Frobenius pair as well. Take $U = \sum_{g \in G} g^2 \in A$.

Let $\phi$ be the natural inclusion from $A$ to $B$.
Let  $\star: B \to B$ be the antipode map, sending~$g$ to~$g^{-1}$.
This map preserves the center, so we have a map $\star: A \to A$,
compatible with the inclusion $\phi$.

\begin{theorem}\label{HG}
The data above form a semi-simple equipped Cardy--Frobenius algebra
over ${\mathbb K}$.
\end{theorem}
\begin{proof}

First let us show that
$U^2 = K_A^\star$.
For a conjugation class $\alpha \subset G$ let
$E_\alpha = \sum_{g \in \alpha} g$. Then~$E_\alpha$ form a basis of $A$.

Note that
$(E_\alpha,E_{\alpha^{-1}})_A = |\alpha|$
and $(E_\alpha, E_\beta)_A = 0$ otherwise, hence as $E_\alpha^\star = E_{\alpha^{-1}}$ we have
\begin{gather*}
K_A^\star = \frac{1}{|G|}\sum_\alpha \frac{E_\alpha^2}{|\alpha|} =
\frac{1}{|G|}\sum_\alpha \sum_{g,g' \in \alpha} \frac{g g'}{|\alpha|} \\
\phantom{K_A^\star}{}  = \frac{1}{|G|}\sum_\alpha \sum_{g \in \alpha, \, h \in G} \frac{g \, h^{-1} gh}{|\alpha|}\frac{|\alpha|}{|G|} =
\sum_{g,h \in G} g h^{-1} gh = \sum_{a,b \in G} a^2 b^2 = U^2,
\end{gather*}
where $a = gh^{-1}$, $b = h$.

Also we have $K_B^\star = \sum_{x,y \in G} l_B(x y^{-1}) xy =  \sum_{g \in G} g^2 = U$.
It remains to prove the Cardy condition  $l_A(\phi^*(x)\phi^*(y)) =\tr W^B_{x,y}$.

We have $\tr W^B_{x,y} = |\{g | xgy = g\}| = |\{g | y = g^{-1}x^{-1} g\}|$.
This number is zero if $x^{-1}$ and $y$ are in dif\/ferent conjugation classes,
and it is $\frac{|G|}{|\gamma|}$ if $x^{-1}$  belongs to the conjugation class $\gamma$ of $y$.
On the other hand, $\phi^*(y) = \sum_{g \in G} g^{-1} y g  =
\frac{|G|}{|\gamma|}\sum_{h \in \gamma} h = \frac{|G|}{|\gamma|}
E_\gamma$.
Thus the number $l_A(\phi^*(x)\phi^*(y))$  is exactly the same.
\end{proof}

We denote the constructed algebra by  $H_{\mathbb{K}}^G$ and call it {\em the regular algebra of $G$}.
This algebra is semi-simple due to semi-simplicity of the group algebra.
Our next aim is full description of~$H^G_\mathbb{R}$ and
$H^G_\mathbb{C}$.

\subsection[Classification of complex semi-simple equipped Cardy-Frobenius algebras]{Classif\/ication of complex semi-simple equipped\\ Cardy--Frobenius algebras}

We call a complex Cardy--Frobenius
algebra $((A,l_A),(B,l_B),\phi)$
\textit{pseudoreal} if $A=A_R\otimes\mathbb{C}$, $B=B_R\otimes\mathbb{C}$ and $\phi=\phi_R\otimes\mathbb{C}$ where $A_R$, $B_R$ are real algebras and $\phi_R:A_R\rightarrow B_R$ is an homomorphism.
It appears (see Theorem~\ref{thc}) that any
simple  complex equipped  Cardy--Frobenius
algebra is pseudoreal.

Let $\mathbb{D}$ be a division algebra over $\R$, that is, $\R$, $\C$ or
$\H$. For each $\mathbb{D}$ introduce a family of real semi-simple equipped Cardy--Frobenius algebras.

Namely, let $n$ be an integer, $\mu \in \C\setminus 0$, put $d = \dim_\R \mathbb{D}$.
Introduce{\samepage
\begin{gather*}
B_R = \mat_n(\mathbb{D}),  \qquad l_{B_R}(x)  = \mu  \, \re \tr (x)  \quad  \mbox{for}\ \ x \in \mat_n(\mathbb{D}),\\
A_R = Z(\mathbb{D}),  \qquad    l_{A_R}(a)  = \mu^2 \re (a)/d      \quad \mbox{for}\  \ a \in Z(\mathbb{D}), \qquad \phi_R(a) = a \,\id \in Z(B_R).
\end{gather*}
For} $z \in \mathbb{D}$ by $\overline{z}$ denote the conjugated element.
The involution $\star_R$ is def\/ined by $a^{\star_R} = \overline{a}$
for $a \in A_R$,
$x^{\star_R} = \overline{x}^t$
for $x \in B_R$, where ${}^t$ means transposition of a matrix.
Now take
\begin{equation}\label{defu}
U_R = \frac{2-d}{\mu} \in A_R.
\end{equation}

Denote this data $((A_R,l_{A_R}),(B_R,l_{B_R}),\phi_R,U_R,\star_R)$ by $H^\mathbb{D}_{n,\mu}$.

\begin{proposition}
The object $H^\mathbb{D}_{n,\mu}$ is a semi-simple real equipped Cardy--Frobenius
algebra.
\end{proposition}

\begin{proof}
Introduce a natural projection $Z: {\mathbb D} \to Z({\mathbb D})$ sending
$x \in \mathbb D$ to $x$ for
${\mathbb D}  = \R,\,\C$ and to $\re (x)$ for ${\mathbb D} =\H$.
We have $\phi^*_R (x) = Z(\tr(x)) d /\mu$. On the other
hand,
a direct calculation in the standard basis shows that
\[
\tr W^B_{x,y} = \tr W^{\mathbb D}_{\tr(x), \tr (y)} = d \re\left(
Z(\tr(x)) Z(\tr(y))\right),
\]
so $(\phi^* (x),
\phi^* (y)) = \tr W^B_{x,y}$.

Another observation is that
 $K_{B_R}^\star =K_\mathbb{D}  \id/\mu$, where  $K_\mathbb{D}$ is the
Casimir element of $\mathbb D$ with respect to the form $(a,b) =
a\overline{b}$. We have $K_\mathbb{D} = (2-d)$, so $\phi(U_R) = K_B$.
At last, $K_{A_R} = d/\mu^2$ for $A_R=\R$ and $K_{A_R} = 0$ for $A_R=\C$, so $U_R^2 = K_{A_R}$.
\end{proof}

\begin{theorem}\label{thc} Any semi-simple equipped
Cardy--Frobenius algebra
$((A,l_A), (B,l_B), \phi, U, \star)$
over~$\mathbb{C}$ is a direct sum of $H^{\mathbb{D}_i}_{n_i,\mu_i} \otimes \C$ and
${\rm Ker} (\phi)$.
\end{theorem}

\begin{proof}
We deduce the statement from a classif\/ication
of  equipped Cardy--Frobenius algebras performed in~\cite{AN}.
To identify $H^{\mathbb{D}_i}_{n_i,\mu_i} \otimes \C$ with the algebras
introduced there
let us describe  $H^{\mathbb{D}_i}_{n_i,\mu_i} \otimes \C$  in detail.
\begin{itemize}\itemsep=0pt
\item If $\mathbb{D}=\R$, then $A \cong \mathbb{C}$  with the identical
involution $x^\star = x$ and the linear form $l_{A} (z)=\mu^2 z$,
$U=\frac{1}{\mu}\in A$;
$B \cong \mat_n(n,\C)$ equipped with
the involutive anti-automorphism $\star: X\mapsto X^t$,
and the linear form $l_{B} (X)=\mu \tr X$.
The homomorphism $\phi:A\to B$ sends the unit to
the identity matrix;
\item If $\mathbb{D} = \C$, then $A \cong \mathbb{C}\oplus\mathbb{C}$ with the
involution $(x,y)^\star=(y,x)$ for $(x,y)\in \mathbb{C}\oplus\mathbb{C}$
and the linear form  $l_A(x,y)=\mu^2 (x+y)/4$,  $U=0$;
$B \cong \mat_n(n,\C) \oplus\mat_n(n,\C)$ with the linear
form $l_B(X,Y)=\mu (\tr X + \tr Y)/2$ and the involutive anti-automorphism
$\star: (X,Y)\mapsto (Y^t,X^t)$.
The  homomorphism $\phi:A\to B$ is given by the equality
$\phi(x,y)=(xE,yE)$.
\item If $\mathbb{D}=\H$, then $(A,l_A,\star)$
is the same as for $\mathbb{D}=\R$, but
$U=-\frac{2}{\mu}\in A$; $B\cong \mat_{2m}(\C)$
with the linear form $l_B(X)=\mu \tr X/2$. A matrix $X\in B$ we
can present in block form as $X=
\left(%
\begin{array}{cc}
  m_{11}&  m_{12}
   \\
  m_{21}& m_{22} \\
\end{array}%
\right)$.
Then the involute anti-automorphism $\star: X\mapsto X^\tau$
is given by the formula
 $X^\tau=
\left(%
\begin{array}{cc}
  m_{22}^t & -m_{12}^t
   \\
 -m_{21}^t & m_{11}^t \\
\end{array}%
\right),$
in other words,  $X^\tau$ is the matrix adjoint to $X$ with respect
to a natural symplectic form.
The homomorphism $\phi:A\to B$ sends the unit to
the identity matrix.\hfill \qed
\end{itemize}\renewcommand{\qed}{}
\end{proof}

\begin{note}
There are real equipped Cardy--Frobenius algebras not isomorphic to
$H^{\mathbb{D}}_{n,\mu}$. A simplest example is a generalization of
$H^{\mathbb{R}}_{n,\mu}$ with the same Cardy--Frobenius structure,
but the involu\-tion~$\star$ def\/ined as the transposition with respect to
a bi-linear form with non-trivial signature.
\end{note}

\subsection{Description of complex regular algebra}

Let us denote complex representations of $G$ by capital Latin letters (such as
$V$) and real representation by Greek letters (such as $\pi$).
Note that irreducible real representations can be split into 3 types
 (see \cite{FH}):
\begin{itemize}\itemsep=0pt
\item Real type:  $\ens (\pi) = \R$ and $\pi \otimes_\R \C$ is irreducible;
\item Complex type: $\ens (\pi) = \C$, $\pi \otimes_\R \C \cong V^+ \oplus V^-$, $V^+$ is not isomorphic to  $V^-$ (but  $V^+ \cong
\left( V^-\right)^*$);
\item Quaternionic  type: $\ens (\pi) = \H$, $\pi \otimes_\R \C \cong V^0 \oplus V^0$.
\end{itemize}

Let ${\rm Ir}_\mathbb{D}(G)$ be the set of isomorphism classes of corresponding
irreducible real representations.
\begin{theorem}\label{dec}
We have isomorphisms of equipped algebras
\begin{gather}\label{eqr}
H_G^\mathbb{R} \cong \bigoplus_{\mathbb{D}=\mathbb{R},\mathbb{C},\mathbb{H}}
\bigoplus_{\pi\in {\rm Ir}_\mathbb{D}(G)}
H^\mathbb{D}_{\frac{\dim\pi}{\dim {\mathbb D}},\frac{\dim\pi}{|G|}},
\\
\label{eqcc}
H_G^\mathbb{C} \cong \bigoplus_{\mathbb{D}=\mathbb{R},\mathbb{C},\mathbb{H}}
\bigoplus_{\pi\in {\rm Ir}_\mathbb{D}(G)}
H^\mathbb{D}_{\frac{\dim\pi}{\dim {\mathbb D}},\frac{\dim\pi}{|G|}} \otimes \mathbb{C}.
\end{gather}
\end{theorem}

\begin{proof} Let us show~\eqref{eqr}, then~\eqref{eqcc} follows
because
$\C[G] \cong \R[G] \otimes \C$.
 By the Wedderburn theorem
\begin{gather}\label{wed}
\R[G] \cong \bigoplus_{\mathbb D}
\bigoplus_{\pi \in {\rm Ir}_\mathbb{D}(G)} \mat_{\frac{\dim\pi}{\dim\mathbb{D}}}(\mathbb{D}),
\end{gather}
where the map $\R[G] \to  \mat_{\frac{\dim\pi}{\dim\mathbb{D}}} (\mathbb{D})$ is the action of $\R[G]$ on
$\pi \cong \mathbb{D}^{\frac{\dim\pi}{\dim\mathbb{D}}}$.  Due to the classif\/ication theorem it is enough
to identify the map $\star$ and the constant $\mu$ on each summand with the same in
$H^\mathbb{D}_{n, \mu}$.

Concerning $\star$, choose an invariant scalar product on $\pi$. As $\pi$ is irreducible,
this invariant bilinear form is unique up to a scalar. Then $\star$ is just the conjugation
with respect to this form. By the orthogonalization process we can suppose that
$\pi \cong\mathbb{D}e_1 \oplus \dots \oplus\mathbb{D} e_m$, where $\{e_l\}$ is a  set of orthogonal vectors.

Note that  $\mat_{\frac{\dim \pi}{\dim\mathbb{D}}}(\mathbb{D})$ is the tensor product
of its subalgebras  $\mat_{\frac{\dim\pi}{\dim\mathbb{D}}}(\mathbb{R})$ and $\mathbb{D}$. In the basis $\{e_l\}$
we identify the action of $\star$ on  $\mat_{\frac{\dim\pi}{\dim\mathbb{D}}}(\R)$ with the matrix transposition.
For $\mathbb{D}$ note that it also acts by right multiplication, an this action commutes with the action of
$\R[G]$, hence this right action preserves the bilinear form up to a scalar. Then it
follows that the set $e_l$, $ie_l$ (for $\mathbb{D}\supset \C$), $je_l$ and $ke_l$ (for $\mathbb{D}=\H$)
form an orthogonal basis. In this basis imaginary elements of $\mathbb{D}$ act by skew-symmetric matrices,
so $\star$ acts on $\mathbb{D}$ as the standard conjugation.

It remains to f\/ind $\mu$ for a summand corresponding to each irreducible real representation~$\pi$.
Let $e_\pi\in A$ be the idempotent corresponding to $\pi$. It acts on
$\R[G]$ by projection onto the corresponding summand in~\eqref{wed}. Then
from the def\/inition
of the regular algebra we have $l_B(\phi(e_\pi)) = {\rm Tr}_{\R[G]} e_\pi /|G|
= \frac{(\dim \pi)^2}{|G|\dim {\mathbb D}}$.
But from the def\/inition of $H^{\mathbb D}_{n,\mu}$ we have
$l_B(\phi(e_\pi)) = \mu n = \mu  \dim \pi/ \dim {\mathbb D}$.
So $\mu = \frac{\dim \pi}{|G|}$.
\end{proof}

\begin{corollary}[cf.~\cite{FH}]\label{SFi}
Let $\pi \in {\rm Ir}_{\mathbb D}(G)$ be an irreducible real representation of $G$. Then $\tr (U)$
on $\pi$
is equal to $(2-\dim {\mathbb D})|G|$
\end{corollary}

\begin{proof}
The element $U$ acts on $\pi \in {\it Ir}_{\mathbb D}(G)$ by the same scalar as on
$H^{\mathbb D}_{\frac{\dim\pi}{\dim {\mathbb D}},\frac{\dim\pi}{|G|}}$.
Substituting the def\/inition~\eqref{defu} for $U$
and multiplying by $\dim \pi$,
we obtain the proposed formula.
\end{proof}

Such an element $U$ is known as  {\em Frobenius--Schur indicator} (see~\cite{FH}). It provides
an easy way to determine type of $\pi$.
It appears in a similar situation in conformal f\/ield theories
(see \cite{B, FFFS}).

\begin{note}
Note that Corollary~\ref{SFi} is applicable to a complex representation in
the same way. Indeed any irreducible complex representation $V$ can be
obtained as a summand in $\pi \otimes \C$ for a real irreducible
representation $\pi$. Then the action of $U$ on $V$ also determines type of $\pi$.
\end{note}

\section[Cardy-Frobenius algebras of representations]{Cardy--Frobenius algebras of representations}\label{s4}

\subsection[Cardy-Frobenius algebra of a representation]{Cardy--Frobenius algebra of a representation}\label{s4.1}

By $\F$ denote a f\/ield.
Let $\rho$ be a f\/inite-dimensional representation (possibly reducible) of a~f\/inite group $G$ over $\F$,
suppose that $|G|$ is not divisible by ${\rm char}\, \F$. Fix an invertible element $\Mu \in  Z(\F[G])$.
Set $A = Z(\F[G])$ with $l_A$ sending an element $x$ to the value of
$\Mu^2 x/|G|$ at the unit element.
Let $B = \ens_G (\rho)$ be the
algebra of intertwining operators on $\rho$ with $l_B (x) = {\rm Tr}_\rho \left(\Mu x/|G|\right)$.
As the center of $\C[G]$ acts on $\rho$ by intertwining operators, we have a natural map $\phi: A \to B$.

\begin{proposition}
The data above form a semi-simple complex Cardy--Frobenius algebra. Let us denote it by $H^\rho_\Mu$.
\end{proposition}

\begin{proof}
As $\Mu$ is invertible, the scalar product on $A$ and $B$ def\/ined by $l_A$ and $l_B$
are non-degenerate, so we have two Frobenius algebras.

Note that $\F$ can be included to an algebraically closed f\/ield $\oF$ and that is  enough
to check the Cardy condition over $\oF$.
The algebra $\oF[G]$ is generated by orthogonal idempotents
$\{e_i\}$,
 corresponding to  irreducible complex representations $\rho_i$ over $\oF$.

The element $\Mu$ acts on $\rho_i$ by a scalar,
let us denote it by $\Mu(\rho_i)$.
Our aim is to show  that the Cardy condition is equivalent to
\begin{gather}\label{CCh}
(e_i, e_i)_A = \Mu^2(\rho_i) \left(\frac{\dim \rho_i}{|G|}\right)^2,
\end{gather}
then we know it for the regular representation, so it
follows for the general case.

If $\rho = \sum \rho_i^{\oplus m_i}$ then $B \cong
\oplus_{i=1}^s\mat_{n_i}(\oF)$.
Note that for $x = \sum x_i \in \oplus_i \mat_{n_i}(\oF)$
and a similar element $y = \sum y_i$
we have
 $\tr W^B_{x,y} =\sum_i \tr x_i \tr y_i$.

Note that $\Mu e_i = \Mu(\rho_i) e_i$, therefore we have
\[
\phi^*(x) = |G| \sum_{i=1}^s e_i \frac{\tr x_i}{\Mu (\rho_i) \dim \rho_i},
\]
so $(\phi^*(x), \phi^*(y))_A = |G|^2 \sum_i (e_i,e_i) (\dim \rho_i)^2  \tr x_i \tr y_i$.
Comparing it with $\tr W^B_{x,y}$ for arbitrary~$x$ and~$y$, we deduce the condition~\eqref{CCh}.
\end{proof}

\subsection[Equipped Cardy-Frobenius algebra of a real representation]{Equipped Cardy--Frobenius algebra of a real representation}\label{s4.2}

Now suppose that $\F = \R$, so $\rho$  is a real representation of $G$.
Then
there is a non-degenerate
symmetric invariant bilinear form on $\rho$.
So for any
operator $x \in \ens (\rho)$
there exists a unique adjoint operator $x^\tau \in \ens (\rho)$.
The map sending $x$ to $x^\tau$ is an anti-involution of $\ens (\rho)$,
preserving the subalgebra  $\ens_G (\rho)$. Thus we obtain a
map $\star : \ens_G(\rho) \to \ens_G(\rho)$.

As before, the involution on $A=\C[G]$ is def\/ined by sending
$g \to g^{-1}$, and $U = \sum_{g \in G} g^2 \in A$. Also we suppose that
$\Mu$ is invariant with respect to this involution, that is,
$\Mu(g) = \Mu(g^{-1})$. Note that such $\Mu$ form a vector space
with dimension equal to the number of real representations of~$G$.

\begin{theorem}\label{rct}
The data above makes $H^\rho_\Mu$  a semi-simple complex  equipped Cardy--Frobenius
algebra.  Moreover, we have  $\rho \cong \bigoplus_{\pi \in {\rm Ir}(G)} n_{\pi} \pi$
and
\[
H^\rho_\Mu \otimes \C \cong \bigoplus_{\mathbb{D }= \R,\,\C,\,\H} \bigoplus_{\pi \in {\rm Ir}_\mathbb{D}(G)}
H^\mathbb{D}_{n_{\pi}, \Mu(\pi) \frac{\dim \pi}{|G|}} \otimes \C.
\]
\end{theorem}

\begin{proof}
The involution $\star$ on $A$ is compatible with
the involution $\star$ on $B$ because sending $g \to g^{-1}$ corresponds to the action on
the dual representation, and this action can be expressed by adjoint
operators with respect to an invariant bilinear form.

By Maschke's theorem we have $\rho \cong \bigoplus_{\pi \in {\rm Ir}(G)} n_{\pi} \pi$.
Note that there is a unique up to a scalar symmetric invariant bilinear form
on a real irreducible representation. Therefore,
this sum can be written as an orthogonal sum such that the restriction of
the form to
each irreducible representation coincides with the scalar product up to a sign.
Then it follows similarly to the case of regular representation that
$H^\rho_\Mu \otimes \C \cong   \bigoplus H^\mathbb{D}_{n_{\pi}, \Mu(\pi) \frac{\dim \pi}{|G|}} \otimes \C$,
in particular, $\left(K_A^\star\right)^2 = U = K_B^\star$ for $H^\rho_\Mu \otimes \C$ and therefore for
$H^\rho_\Mu$ itself.
\end{proof}

\begin{corollary}
Any semi-simple equipped complex Cardy--Frobenius algebra can be constructed from a real representation
of a finite group.
\end{corollary}

\begin{proof}
By $Q_8$ denote the group of 8 elements $\pm 1$, $\pm i$, $\pm j$, $\pm k$ with the natural quaternion multiplication.
It has one 4-dimensional quaternionic type representations over $\R$ and four 1-dimensional real type representations.
Then the group $Q_8 \times \Z/3\Z$ has all the types of representations, so using Cartesian product of several copies
of this group by Theorem~\ref{rct} and Theorem~\ref{dec} we obtain an
arbitrary semi-simple equipped complex Cardy--Frobenius algebra as a
summand.
\end{proof}

\subsection{Group action case}
A particular case of this construction was already discovered in
\cite{AN3}.
Suppose that the group $G$ acts on a f\/inite set $X$. Let $\pi_X = {\mathbb R} X $ be
the real representation of $G$ in the vector space formed by formal
linear combinations of the elements of $X$.

Let $H^{\pi_X}_e = ((A,l_A),(B,l_B),\phi,U,\star)$. Then
an explicit construction of $B$ is proposed in \cite{AN3}.

The group $G$ acts on $X^n=X\times \dots \times X$ by formula
$g(x_1,\dots,x_n)= (g(x_1),\dots,g(x_n))$.  Let
$\mathcal{B}_n=X^n/G$. By $\Aut\bar x$
denote the
stabilizer of element $\bar x \in X_n$.
Indeed for $\bar x=(x_1,\dots,x_n)$ we have
$\Aut\bar x =\cap_i\Aut x_i$.
Cardinality of this subgroup $|\Aut\bar x|$ depends only on
the orbit of $\bar x$, so we consider it as a function on
$\mathcal{B}_n$.

By $B_X$ denote  the vector space generated by $\mathcal{B}_2$. The
involution $(x_1,x_2)\mapsto(x_2,x_1)$ generates the involution
$\star_X :B_X\rightarrow B_X$. Introduce a bi-linear and a three-linear
form on $B_X$ as follows:
\[
(b_1,b_2)_X = \frac{\delta_{b_1,b_2^\star}}{|\Aut b_1|},
\qquad
(b_1,b_2,b_3)_X = \sum_{(x_1,x_2) \in b_1,\ (x_2,x_3) \in b_2, \ (x_3,x_1) \in b_3}
\frac{1}{|\Aut (x_1,x_2,x_3)|}.
\]
Def\/ine a multiplication on $B_X$ by $(b_1b_2,b_3)_X=
(b_1,b_2,b_3)_X$. The element $e = \sum_{x\in X} (x,x)$ is a~unit of $B_X$.
At last, let  $l_{B_X}(b)=(b,e)_X$.

\begin{theorem}
We have an isomorphism $B \cong B_X$ identifying  $l_B$ with  $l_{B_X}$ and
$\star$ with $\star_X$.
\end{theorem}

\begin{proof}
Essentially, it was done in \cite{AN3}. Elements $(x_1,x_2) \in X\times X$
enumerates matrix units $E_{x_1,x_2} \in \ens(\pi_X)$, so to any orbit
$b\in\mathcal{B}_2$
we correspond the operator $\sum_{(x_1,x_2) \in b} E_{x_1,x_2} \in \ens_G(\pi_X)$.
One can check by a direct computation that this map is an algebra homomorphism and that
the trace $l_B$ can be written as  $l_{B_X}$.
Clearly, this homomorphism is injective and, as the operator $\sum_{g \in G} g$ on  $\ens(\pi_X)$
 is the projection to the subspace of invariants, it is surjective,
so we have $B_X \cong B$.
At last, the involution $\star_X$ corresponds to transposition of a matrix in
the natural orthonormal basis of $\pi_X$, hence it corresponds to $\star$.
\end{proof}

\begin{note}
This construction def\/ines a structure of real equipped Cardy--Frobenius algebra on
the Hecke algebra $H\setminus G /H$ for an arbitrary subgroup $H \subset G$.
To this end one can take $X$  to be the left coset
$G/H$ with the
natural action of $G$.
\end{note}

\subsection*{Acknowledgements}

We are grateful to P.~Deligne, B.~Feigin, Yu.~Manin, S.~Shadrin and V.~Turaev for useful discussions.
Part of this work was done during the stays of S.N.\  at
Max-Planck-Institute in Bonn,
he is grateful to MPIM for their hospitality and support.
The work of S.N.\ was partly supported by grants RFBR-11-01-00289,
N.Sh-8462.2010.1 and the Russian government grant 11.G34.31.0005. The work of S.L.\  was partly supported by grants:
N.Sh-3035.2008.2,  RFBR-09-01-00242, SU-HSE award No.09-09-0009, RFBR-CNRS-07-01-92214, RFBR-IND-08-01-91300,
 RFBR-CNRS-09-01-93106 and P.~Deligne 2004 Balzan prize in mathematics.

\pdfbookmark[1]{References}{ref}
\LastPageEnding

\end{document}